\documentclass[12pt, leqno]{article}
\usepackage{amssymb}
\usepackage{amsmath}
\title{Derived brackets\thanks{This article is the revised version of 
a lecture given at the Euroconference,
    Poisson Geometry, Deformation Quantization and Group
    Representations (PQR 2003) at the Universit{\'e} Libre de Bruxelles,
    18-22 June 2003.}}

\author{Yvette Kosmann-Schwarzbach\\
Centre de Math{\'e}matiques,
U.M.R. 7640 du CNRS,\\
{\'E}cole Polytechnique\\
F-91128 Palaiseau, France\\
{\tt yks@math.polytechnique.fr}}

\date{}

\def\T*{T^*M}
\def\t{\tilde}

\def\d{\rm d}

\def\End{\rm End}
\def\W{\Omega^{\bullet}(M)}

\def\o{\otimes}
\def\M{\mathcal M}
\def\bra{[~,~]}
\def\C{C^{\infty}(M)}

\newtheorem{theorem}{Theorem}

\newtheorem{definition}{Definition}
\newtheorem{proposition}{Proposition}
\newtheorem{corollary}{Corollary}

\numberwithin{equation}{section}
\numberwithin{theorem}{section}
\numberwithin{definition}{section}
\numberwithin{proposition}{section}
\numberwithin{corollary}{section}

\begin{document}

\maketitle

\abstract{\it We survey the many instances of derived bracket
  construction in differential geometry, Lie algebroid and Courant
  algebroid theories, and their properties.
We recall and compare the constructions of Buttin and of 
Vinogradov, and we prove that the Vinogradov bracket is the
  skew-symmetrization of a derived bracket.
Odd (resp., even) Poisson brackets on supermanifolds are derived
  brackets of canonical even (resp., odd) Poisson brackets on their 
cotangent bundle (resp., parity-reversed cotangent bundle). 
Lie algebras have analogous properties, and the theory of Lie algebroids
  unifies the results valid for manifolds on the one hand, and for Lie
  algebras on the other.
We outline the role of derived brackets in the
  theory of ``Poisson structures with background''.}

\section*{Introduction}
On any graded differential Lie algebra, $(\mathfrak A, $\bra$, D)$, with
bracket of degree $n$,
one can consider the bilinear map,
$$
(a,b) \in \mathfrak A \times \mathfrak A \mapsto (-1)^{n+|a|+1} [Da,b] 
\in \mathfrak A \ ,
$$
where $|a|$ is the degree of $a$.
This is what is called the {\it derived bracket} of $\bra$ by
$D$. (See \cite{yks96}.)
It is not in general a graded Lie bracket because it is
not, in general, skew-symmetric. However
it does satisfy the Jacobi identity in the form \eqref{jac} below,
therefore $\mathfrak A$ with a derived bracket 
is a graded version of what Loday calls a {\it Leibniz
algebra}, which we prefer to call a {\it Loday algebra}. 
Since the derivation $D$ is odd, passing from the original 
bracket to the derived
bracket
turns an even (resp., odd) Lie bracket into an odd (resp., even)
Loday bracket.

In applications, $D$ will most often be the interior
derivation by an odd element, $d \in \mathfrak A$, of square $0$,
$$
[d,d]=0 \ .
$$
The derived bracket can then be written simply as 
$$
(a,b) \in \mathfrak A \times \mathfrak A \mapsto [[a,d],b]
\in \mathfrak A \ .
$$
Whenever $(\mathfrak A, \bra)$ is
abelian, the derived bracket
is a genuine graded Lie bracket, 
and this property is essential for the applications that we shall describe.

Instances of derived brackets, though they do not bear that name, 
appeared in various contexts:
in an article on formal non-commutative geometry by 
Gel'fand, Daletskii and Tsygan \cite{GDT},
in papers by physicists
on the BRST quantization, which are too numerous to be exhaustively
cited here, see
\cite{BM},
in the 1974 article by Buttin \cite{B},
in early work of A. M. Vinogradov \cite{V}
who introduced a very powerful tool under
the unfortunate term of ``lievization'', in unpublished papers of Ted
Voronov, and certainly in other sources of which I am not aware.
The notion was formalized in unpublished notes of
Koszul dating from 1990 \cite{K}, which he communicated to me in 1994. 
There followed the article \cite{yks96}
where I placed Koszul's construction in the framework of graded Loday
algebras, proved the main properties of the general construction,
and gave examples from Poisson geometry and Lie bialgebra theory. 
Related results were obtained independently by Daletskii and
Kushnirevitch \cite{DK}. Articles \cite{yksgoslar} and \cite{yksoddeven}
contain a summary of results and describe applications to gauge
Lie algebras in various field theories.

After briefly reviewing the general notion, I shall try to describe 
enough old and new examples 
of derived brackets to convince the reader of their ubiquity and importance.
The following discussion may seem very formal, but the
general results on the derived brackets of Lie brackets
briefly recalled in Section~\ref{derived} of 
this survey, and more
generally of Loday brackets, are powerful tools for proving
non-trivial properties of brackets and derivations.

\section{Derived brackets}\label{derived}

\subsection{Loday brackets}
Loday algebras were introduced (in the ungraded case) by Jean-Louis
Loday under
the name {\it Leibniz algebras} \cite{L}.
We define a {\it Loday algebra of degree} $n$ 
as a graded vector space $V$ over a field $R$
of characterisitic $\neq 2$
(or just a module over a commutative ring), 
equipped with an $R$-bilinear map, $[~,~]: V \otimes V \to V$, satisfying
the Jacobi identity in the form,
\begin{equation}\label{jac}
[a,[b,c]] = [[a,b],c] + (-1)^{(n + |a|)(n + |b|)}[b, [a,c]] \ ,
\end{equation}
for all $a, b$ and $c \in V$, where $|a|$ denotes the degree of $a \in V$. 
Whenever the bracket $[~,~]$ is graded skew-symmetric, a 
Loday algebra is just a graded Lie algebra.
In what follows, we shall often omit the word ``graded''.

\subsection{Definition of derived brackets}\label{def}
In \cite{yks96} 
(see also \cite{yksgoslar} 
\cite{yksoddeven}), we defined a general notion of {\it derived
  brackets} of Loday brackets, and we
proved some simple properties that have far-reaching
consequences. This construction turns an even Loday bracket into an
odd one, and conversely.
Here, for simplicity, we describe this construction in the case
of derived brackets of Lie brackets, 
and, in Theorems \ref{theoremderived} and \ref{abel} below, we 
recall those properties that are most useful in
applications.

\begin{definition} If $(V, [~,~], D)$ is a graded 
differential Lie algebra over $R$ with bracket of degree $n$,
we define the bilinear map $[~,~]_{(D)}: V \otimes V \to V$ by
\begin{equation}
[a,b]_{(D)} = (-1)^{n+|a|+1} [Da,b] \ ,
\end{equation}
for $a$ and $b \in V$,
and we call it 
the {\em derived bracket} of $[~,~]$ by $D$.
\end{definition}

\begin{theorem}\label{theoremderived}
(i) The derived bracket of a Lie bracket of degree $n$ is a Loday bracket
of degree $n+1$.

(ii)
The map $D$ is a morphism of Loday algebras from $(V, \bra_{(D)})$ to
$(V,\bra)$.

(iii)
The map $D$ is a derivation of the Loday bracket $\bra_{(D)}$.

(iv) The restriction of the derived bracket to any {\em abelian
  subalgebra}, $V_0$, 
of $(V, [~,~])$, such that $[DV_0,V_0]\subset V_0$,  
is a Lie bracket of degree $n+1$. 

(v) The bracket $[~,~]_{(D)}$ induces a Lie bracket of degree $n+1$
on the quotient space of $V$ by the image of $V$ under $D$, $V/D(V)$.
\end{theorem}

More generally, we can consider the derived bracket by any derivation
of odd degree and of square $0$.

\medskip

\noindent{\bf Example} 
The problem of extending the Poisson bracket of functions into an even
bracket on the algebra of all differential forms was a long-standing
problem until the mid-nineties. 
In \cite{yks96}, we proved 
that such an extension can be easily defined, but only as a
Loday bracket. Let $P$ be a Poisson bivector on a smooth manifold
$M$, and let $\bra^P$
be the Koszul bracket of differential forms \cite{Ko}. 
The derived bracket of $\bra^P$ by the de Rham differential is an even
Loday bracket 
which extends the Poisson bracket of functions. 
We observed that one of the brackets defined on symplectic
manifolds by Michor in \cite{M},
denoted there by $\{~,~\}^2$, coincides with 
this derived bracket.
Shortly after that,
Grabowski 
proved that an extension as
an even graded Lie bracket can be defined, but this bracket is not a
biderivation of the algebra \cite{grabowski}.

\subsection{The case of an interior derivation}
In many applications, the derivation $D$ is an interior derivation of
$(V, \bra)$, $a \mapsto~[d,a]$, where $d$ is an element of square $0$
in $(V, \bra)$.

\medskip

\noindent{\it Notation} 
If $D$ is the interior derivation by an element $d$ in the Lie algebra
$V$, we denote the corresponding derived bracket 
simply by $[~,~]_d$.

\begin{theorem}\label{abel}
If $D$ is the interior derivation of $(V,\bra)$
by an element $d \in V$ such that $|d|+n$ is odd and $[d,d]=0$, 
the derived bracket is 
\begin{equation}\label{interior}
[a,b]_d = [[a,d],b] \ ,
\end{equation}
for $a$ and $b \in V$. Both $a \mapsto [d,a]$ and $a \mapsto [a,d]$
are morphisms from $(V, \bra_d)$ to $(V,\bra)$. 
\end{theorem}

The proof is obtained by a simple calculation.

\medskip

\noindent{\bf Example} Let $P$ and $\bra^P$ be as in the Example of
Section \ref{def}. Assume that $P$ is non-degenerate, 
with inverse symplectic form $\omega$. Then 
the de Rham differential is the interior derivation $[\omega, \cdot~]^P$,
as was proved in, {\it e.g.}, \cite{yksm}. Therefore, in this case,
the derived bracket of forms
$\alpha$ and $\beta$ is
 equal to $[[\alpha, \omega]^P,\beta]^P$.

\subsection{Skew-symmetrization of derived brackets}
The skew-symmetrization of the derived bracket $\bra_{(D)}$, 
which we denote by
$\bra_{(D)}^{-}$, can be expressed as
\begin{equation}\label{skew}
[a,b]_{(D)}^{-}= \frac{1}{2}\left([a,Db]-(-1)^{n+|a|}[Da,b]\right) \ ,
\end{equation}
while, in the case of an interior derivation, it satisfies
\begin{equation}\label{skewint}
[a,b]^{-}_d= \frac{1}{2}\left([[a,d],b] - (-1)^{|b|} [a, [b,d]]\right) \ .
\end{equation}
In general, this bracket, obtained by skew-symmetrizing a Loday
bracket, no longer satisfies the Jacobi identity 
\eqref{jac}: a defect in the
Jacobi identity appears, so the skew-symmetrized bracket is not a Lie
bracket. We shall give examples below.
 
\subsection{Notations}
When dealing with the geometric objects defined on a smooth manifold
or supermanifold, we shall assume all fields to be smooth, and we shall
often abbreviate {\it vector field} and
{\it multivector field} to {\it vector} and {\it multivector}, 
respectively, and
{\it differential form} to {\it form}.
If $M$ is a manifold, 
we denote the exterior algebra of multivector fields by
$V^{\bullet}(M)   = \bigoplus_{p \geq 0} V^{p}(M)
= \Gamma(
\bigoplus_{p\geq 0} \bigwedge^p(TM))$, so that $V^1(M)$
is the space of vector fields, and we denote
the exterior algebra of 
differential forms by
$\W= \bigoplus_{q \geq 0} \Omega^{q}(M)   
= \Gamma(\bigoplus_{q \geq 0}  \bigwedge^q(T^*M))$, so that
$\Omega^1 (M)$ is the space of differential $1$-forms.
(If $E$ is a vector bundle over $M$, we denote the space of sections of
$E$ by $\Gamma E$.)
Because the Schouten-Nijenhuis bracket of multivectors \cite{S}
\cite{N55} \cite{N} 
(or Schouten bracket for short) is the prototypical example of a
Gerstenhaber bracket, we use also the term {\it Schouten algebra} to
designate a {\it Gerstenhaber algebra}, {\it i.e.},
an associative, graded commutative 
algebra with an odd Poisson bracket.

The notation $[~,~]$ stands for the graded commutator
of graded endomorphisms of a graded module, unless specified otherwise.

\section{The Cartan formulas and their generalizations}
\subsection{The Cartan formulas}
Everyone is familiar with the formulas to be found in, {\it e.g.}, 
Henri Cartan's celebrated communication presented in 
Brussels half a century ago
\cite{Cartan},
\begin{equation}\label{cartaneq}
[{\d},{\d}]=0, 
~~~[i_x,i_y]=0,~~~L_x=[i_x,{\d}], ~~~[L_x,{\d}]=0,~~~ [L_x,i_y]=i_{[x,y]} \ .
\end{equation}
Here $x$ and $y$ are vector fields on a manifold, ${\d}$ is the de Rham
differential,
$i _x$ is the interior product by $x$ and $L_x$
is the Lie derivation by $x$, each of these being a derivation of the
algebra of differential forms,
while $[~,~]$ 
is the graded commutator of graded endomorphisms of the
space of differential forms, but $[x,y]$ in the last 
formula denotes the Lie bracket of
vector fields.

\subsection{$(\mathfrak A,\mathcal D)$-structures}
Motivated by the extension of the usual differential calculus to the
calculus of variations (see \cite{GD1980}),
Gel'fand, Daletskii and Tsygan \cite{GDT}, working
with multigraded vector spaces, present the
general theory of an $(\mathfrak A, \mathcal D)$-structure, a graded
Lie algebra $\mathfrak A \oplus \mathfrak g_0 \oplus \mathcal D$,where
$\mathfrak A$ is an abelian Lie algebra, generalizing 
the space of the $i_x$'s, and where $\mathcal D$ is a 
space generated by several commuting elements of square $0$, generalizing 
the one-dimensional vector space generated by the de Rham differential,
while $\mathfrak g_0$ generalizes the Lie algebra of Lie derivations.
For any $d \in \mathcal D$, 
they let 
$$[a,b]_d=[[a,d],b],
$$
for $a, b \in \mathfrak A$, and they
formulate their Theorem 1, stating that this formula defines a 
graded Lie bracket on $\mathfrak A$, and that the map $a \mapsto
[a,d]$ is a morphism of graded Lie algebras, from $\mathfrak A$ with
bracket $[~,~]_d$ to $\mathfrak g_0$ with the original Lie bracket.
It is clear that 
these statements are special cases of the results summarized in 
Theorems \ref{theoremderived} and \ref{abel} above.

In addition, Gel'fand, Daletskii and Tsygan state a
property which is equivalent to the compatibility of 
the brackets $[~,~]_{d_1}$ and $[~,~]_{d_2}$ on $\mathfrak A$, for a pair of
commuting differentials $d_1$ and $d_2$, 
in the sense that the bracket 
$[~,~]_{d_1} + [~,~]_{d_2}$ is a graded Lie bracket on $\mathfrak A$.
Since 
$$
[a,b]_{d_1} +[a,b]_{d_2} =[a,b]_{d_1+d_2} \ ,  
$$
this compatibility follows from the assumptions 
$[d_1,d_1]=[d_2,d_2]=[d_1,d_2]=0$, which imply that $[d_1+d_2,d_1+d_2]=0$.

\subsection{The Cartan formula for multivectors}
The last of the above-mentioned Cartan equations \eqref{cartaneq}, which can
be written  
\begin{equation}\label{cartan}
i_{[x,y]} = [[i_x,{\d}],i_y] \ ,
\end{equation}
expresses the fact that {\it the Lie bracket of vector fields is a derived
bracket}.

It is well-known \cite{T} \cite{O} \cite{Ca} that equation
\eqref{cartan}
is valid more generally for multivectors,
when the bracket on the left-hand side is the Schouten-Nijenhuis
bracket, showing that {\it the Schouten-Nijenhuis bracket of
multivectors is a derived bracket}.

In the sequel of their paper \cite{GDT}, Gel'fand, Daletskii and Tsygan show
that equation \eqref{cartan} 
is valid in the more general case where 
$x$ and $y$ are elements of the 
associative, graded commutative algebra
generated by $\mathfrak A$, and the bracket on the
left-hand side is extended by the bi-derivation property. This 
generalizes the preceding statement concerning multivectors
on manifolds
to the case of 
$({\mathfrak A}, {\mathcal D)}$-structures. In \cite{GDT}, 
further applications are made to the Gerstenhaber algebra 
structure of the
Hochschild cohomology of an associative superalgebra.

\subsection{Some historical comments}
The theory of $({\mathfrak A}, {\mathcal D})$-structures, which was
formulated around 1987,
is closely
related to that of {\it complexes over Lie algebras}, developped
around 1980 by 
Gel'fand and Dorfman \cite{GD1980}, 
as an abstract framework for the variational calculus and the theory of 
integrable systems. In a complex over a Lie algebra, the Cartan
formulas \eqref{cartaneq} are taken as the defining properties.

\medskip

\noindent{\bf Remark} It is easily seen that there is a complex over a
Lie algebra associated to any Lie algebroid, the complex being the
algebra of sections of the exterior algebra of the given vector bundle,
equipped 
with the Lie algebroid differential (see {\it e.g.}, Section \ref{algd}
below). The same conclusion is valid for a Lie-Rinehart algebra.
These are particular cases of complexes over Lie algebras, 
since no associative multiplication is
assumed on the total space of the complex in general. 

\medskip

In a 1986 preprint, Ted Courant and Alan Weinstein 
introduced the notion of {\it Dirac structure on a manifold}, by
imposing an integrability condition on a field of Dirac structures at
each point of the manifold, defined as totally isotropic subspaces in the 
direct sum of the tangent and the cotangent bundle.
Their study was pubished in \cite{CW}.
Later, Courant proved in \cite{C} 
that this integrability condition amounts to a closure condition 
under a bracket on $\Gamma (TM \oplus T^*M)$, this bracket that now
bears his name
being skew-symmetric but not satisfying
the Jacobi identity.
The theory was later developed by Zhang-Jiu Liu, Weinstein and
Ping Xu, who
introduced the more general notion of Courant algebroid, and 
of Dirac structures as integrable subbundles of Courant algebroids \cite{LWX}.
See Section \ref{courantalg} below.
  
Meanwhile, inspired by the finite-dimensional structures first considered by
Cou\-rant and Weinstein in 1986,
Irene Dorfman introduced \cite{Dorf87} \cite{Dorf}
a general notion of Dirac structure in the algebraic framework of 
complexes over Lie
algebras. In both cases the motivation was to unify the pre-symplectic
and Poisson structures (called {\it Hamiltonian structures} in the
infinite-dimensional case).
The equivalence of the definitions in the case of the de Rham complex
over the Lie algebra of vector fields on a smooth manifold is not
explicit in the literature, but is not hard to prove.

\section{The brackets of Buttin, Vinogradov and\\ Courant}

\subsection{The Buttin brackets}
In \cite{B}, an article developped  
by Pierre Molino from the notes left by Claudette Buttin
(1935-1972),
we find a study of the differential operators of all orders 
on the exterior algebra of a module. Buttin calls the order 
the ``type'', to distinguish it from the order in the usual
sense when the exterior algebra is the algebra of forms on a smooth
manifold, and we shall follow this convention.
Recall that an endomorphism of a graded associative algebra is
called a differential operator of type~$0$ if it commutes (in the
graded sense) with the
left multiplication by any element in the algebra, and of type $\leq k$
if its graded commutator with the left multiplication by any element
in the algebra is of type $\leq k - 1$.

Buttin first defines a composition law extending the 
Nijenhuis-Richardson bracket
on the space  $\bigwedge^{\bullet} E ^* \otimes E$ 
of vector-valued forms on a module $E$ \cite{NR} 
to a graded Lie bracket defined on the space of
all {\it multivector-valued forms}, $\bigwedge^{\bullet} E ^* \otimes 
\bigwedge^{\bullet} E$.
This bracket is defined 
by considering the embedding, $i$, of $\bigwedge^{\bullet} E^* \otimes
\bigwedge^{\bullet}E$ 
into the graded Lie algebra of all diffferential operators
on $\bigwedge^{\bullet} E^*$ equipped with the graded commutator.
Let $i_x$ be the interior product of forms by a vector $x \in E$.
For a decomposable multivector, $x = x_1 \wedge \ldots \wedge x_p \in 
\bigwedge ^p E$, we set $i_x = i_{x_1} \ldots i_{x_p}$.
The operator $i_x$ is of type~$p$.

\begin{definition} The {\em embedding} $i$ of 
$\bigwedge^{\bullet} E^* \otimes
\bigwedge^{\bullet}E$ 
into the vector space of all diffferential operators
on $\bigwedge^{\bullet} E^*$ is defined on decomposable elements by
\begin{equation}\label{embed}
i_{\xi \otimes x} (\alpha) = \xi \wedge i_x \alpha \ ,
\end{equation}
for $x \in
\bigwedge^{\bullet} E$, $\xi$ and $\alpha \in \bigwedge^{\bullet} E^*$.
\end{definition}
This embedding restricts to the {\it interior product} of forms by
multivectors on the one hand, and to the {\it exterior product} by forms on
the other hand. We shall use this definition in what follows,
sometimes 
reverting to the notation $e_{\xi}$ instead of
$i_{\xi}$
when $\xi$ is a form acting on forms by exterior product. 
In this notation, $i_{\xi \otimes x} = e_{\xi} \circ i_x$.
The operator $i_{\xi \otimes x}$ is of degree $|\xi| - |x|$, and of type $|x|$.

\medskip

We introduce the notation $\bra^0_B$ for Buttin's algebraic bracket,
in which her definition becomes

\begin{definition}
For $X$ and $Y \in \bigwedge^{\bullet} E^* \otimes
\bigwedge^{\bullet}E$, the bracket $[X,Y]^0_B$ is the element in 
$\bigwedge^{\bullet} E^* \otimes
\bigwedge^{\bullet}E$
such that $i_{[X,Y]^0_B}$ 
is the {\em term of highest type} in $[i_X,i_{Y}]$.
\end{definition}

For $X \in \bigwedge^{q} E^* \otimes 
\bigwedge^{p} E$, $Y \in \bigwedge^{q'} E^* \otimes 
\bigwedge^{p'} E$, the bracket $[X,Y]^0_B$ is in 
$\bigwedge^{q+q'-1} E^* \otimes \bigwedge^{p+p'-1} E$.
We wish to compare this little-known construction with the well-known
notion of the {\it big bracket} \cite{KosSte} \cite{LR} \cite{yksJ}.
The {\it big bracket} is the canonical Poisson structure
on $\bigwedge^{\bullet}(E \oplus E^*)$, 
the even Poisson bracket 
on the cotangent bundle of the odd
supermanifold
obtained from $E$ (or $E^*$) 
by a change of parity (see Section \ref{bigb}). For the properties of
the big bracket, see \cite{yksJ} \cite{BKS}. 
Here and below, we denote it by $\{~,~\}$.
The comparison of Buttin's bracket, $\bra^0_B$, with the big bracket 
is a consequence of the following result.

\begin{theorem}
For any $X$ and $Y \in \bigwedge^{\bullet} E^* \otimes
\bigwedge^{\bullet}E$,
$i_{\{X,Y\}}$ is the term of highest type in $[i_X,i_Y]$.
\end{theorem}

\noindent{\it Proof} 
Both $\{X,Y\}$ and $[i_X,i_Y]$ are $0$ when both arguments $X$ and $Y$
belong to   
$\bigwedge^{\bullet} E$, or to $\bigwedge^{\bullet} E^*$.
If $X \in E$, $Y \in E^*$, then both $i_{\{X,Y\}}$ and $[i_X,i_Y]$ are
multiplication by the scalar $<Y,X>$, obtained from the duality of
$E^*$ with $E$. It is now enough to remark that both expressions are
derivations with respect to $X$ and $Y$. 
On the one hand, we know that the big bracket
satisfies, for $X, Y$ and $Z \in  \bigwedge^{\bullet} E^* \otimes
\bigwedge^{\bullet}E$, 
$$ 
\{X,Y\wedge Z\}=\{X,Y\}\wedge Z + (-1)^{|X||Y|} Y \wedge \{X,Z\} \ .
$$
On the other hand, it follows from the properties of the graded
commutators that
$$
[i_X,i_{Y \wedge Z}]= [i_X,i_Y] \circ i_Z + (-1)^{|X||Y|} i_Y \circ
[i_X,i_Z] \ .
$$
Since the term of highest type of the composition of two
differential operators is the composition of the terms of highest
type, the theorem follows.

\begin{corollary} \label{big}
In all cases, $\bra ^0_B$  
coincides with the big bracket.
\end{corollary}

In particular, if $p=p'=1$, $[i_X,i_Y]$ has only terms of highest type
and therefore 
\begin{equation}\label{bzero}
i_{[X,Y]^0_B}= [i_X,i_Y] \ .
\end{equation}
This property also follows from the fact that the
restriction of the big bracket to the vector-valued
$1$-forms
is equal to the Nijenhuis-Richardson bracket. 
The same
is true of bracket $[~,~]^0_B$ by Corollary \ref{big}.
Therefore, Equation \eqref{bzero} reduces to the 
defining property of the Nijenhuis-Richardson bracket \cite{NR}:
the embedding $i$ maps the Nijenhuis-Richardson bracket on $\bigwedge
^{\bullet}  E^* \otimes E$ to the commutator of operators on 
$\bigwedge^{\bullet}E^*$.

\medskip

\noindent{\bf Remark}
Contrary to a statement in \cite{B}, relation \eqref{bzero} does not hold in
all generality, as can be shown by the calculation of 
$[i_{x\wedge y},
e_{\xi}]$, 
where $x$ and $y$ are vectors, and $\xi$ is a form of degree
$\geq 2$.
In this case 
$$
[i_{x\wedge y},e_{\xi}] = e_{i_y \xi}~ i_x + (-1)^{|\xi|} e_{i_x\xi}
~ i_y - (-1)^{|\xi|} e_{i_{x\wedge y}\xi} \ ,
$$ 
which is the sum of a term of type $1$ and a term of type $0$.
Identifying $X$ with $i_X$, one can write,
$$
[{x\wedge y},{\xi}] = {i_y \xi} \otimes x + (-1)^{|\xi|}  {i_x\xi}
\otimes y
- (-1)^{|\xi|} {i_{x\wedge y}\xi} \ ,
$$
while only the first two terms constitute 
$[x\wedge y , \xi]^0_B$.

Another example is the Buttin algebraic bracket of two
bivector-valued $1$-forms, showing that, in this case also, the result
is the sum of the term of highest type
(which is $3$ in this example) 
and of terms of lower type. Explicitly,
$$
[{\xi_1} \otimes {x_1\wedge y_1}, {\xi_2} \o {x_2\wedge y_2}]
$$
$$
= 
(-1)^{|\xi_2|}~ 
{\xi_1 \wedge {i_{x_1}\xi_2}} \o {y_1 \wedge x_2 \wedge y_2} 
+ 
{\xi_1 \wedge {i_{y_1}\xi_2}} \o {x_1 \wedge x_2 \wedge y_2}
$$
$$
- (-1)^{(|\xi_1|+1)(|\xi_2|+1)} ~
( {\xi_2 \wedge {i_{x_2}\xi_1}} \o {y_2 \wedge x_1 \wedge y_1}
+ (-1)^{|\xi_1|} 
{\xi_2 \wedge {i_{y_2}\xi_1}} \o {x_2 \wedge x_1 \wedge y_1})
$$
$$
+ (-1)^{|\xi_2|+1} 
({\xi_1 \wedge {i_{x_1 \wedge y_1}\xi_2}} \o {x_2 \wedge y_2} +
(-1)^{|\xi_1|(|\xi_2|+1)} 
{\xi_2 \wedge {i_{x_2 \wedge y_2}\xi_1}} \o {x_1 \wedge y_1})
\ .
$$
We observe that the expression of the term of highest type 
coincides with the explicit formula for
the big bracket given in \cite{BKS}. 

\bigskip

Buttin then considers the case of differential operators on 
the exterior algebra,
$\W =
 \Gamma (\oplus_{q \geq 0}\bigwedge^q (T^*M))$, 
of all differential forms on a smooth manifold, $M$. If 
${\d}$ is the de Rham differential, and if $X$ and $Y$ are 
{\it multivector-valued differential forms} on $M$, {\it i.e.},
tensors
skew-symmetric in both their contravariant and their covariant indices,
one can consider
the expression, $[[i_X,{\d}],[i_{Y},{\d}]]$.
She proves 
that there exists 
a differential operator, $\{X,Y\}_B$, 
called the {\it generalized differential concomitant of the
  first kind}, such that  
\begin{equation}\label{concom}
[\{X,Y\}_B,{\d}]  =
[[i_X,{\d}],[i_{Y},{\d}]] \ ,
\end{equation}
which is well-defined when an additional
condition
is imposed on its symbol.
In addition, she shows that, only in certain cases,
there exists a tensor $[X,Y]_B$ such that 
\begin{equation}\label{diffbra}
 [i_{[X,Y]_B},{\d}] = [[i_X,{\d}],[i_{Y},{\d}]] \ .
\end{equation}
These cases are

(0) $X$ is a differential form, and in
this case, $[X,Y]_B$ is not a true differential concommitant 
since it does not involve partial derivatives of the components of
$Y$,

(1) $X$ and $Y$
are multivector fields, and in this case $[~,~]_B$ 
is the Schouten-Nijenhuis bracket as it is now usually defined (differing from
Lichnerowicz's definition \cite{Li} by a sign),

(2) $X$ and $Y$ are vector-valued
differential forms, and in this case $[~,~]_B$ coincides with  
the Fr{\"o}licher-Nijenhuis bracket \cite{FN}, \cite{DVM}.
In fact, setting $L_X=[i_X,\d]$, formula \eqref{diffbra} becomes
\begin{equation}
L_{[X,Y]_B} = [L_X,L_Y] \ ,
\end{equation}
and the Fr{\"o}licher-Nijenhuis bracket is a solution of this
equation. Because, in this case, $i_X, i_Y$ and $\d$ are derivations of
$\Omega^{\bullet}(M)$, we know that this solution is unique.

The brackets of cases (1) and (2), 
each extending the Lie bracket of vector fields,
are thus seen as particular cases of a more general construction. 

\medskip 

\noindent{\bf Remark} 
The embedding $i$ can be considered as an embedding of  
$\bigwedge^{\bullet} E^* \otimes
\bigwedge^{\bullet}E$ 
into the vector space of all diffferential operators
on $\bigwedge^{\bullet} E^* \otimes \bigwedge^{\bullet} E$, 
defined on decomposable elements by
\begin{equation}\label{embedding}
i_{\xi \otimes x} (\eta \otimes y) = 
i_{\xi \otimes x}\eta \otimes y =\xi \wedge i_x \eta \otimes y \ ,
\end{equation}
for $x$ and $y \in
\bigwedge^{\bullet} E$, $\xi$ and $\eta \in \bigwedge^{\bullet} E^*$.
Similarly, there is an embedding $j$ of 
$\bigwedge^{\bullet} E \otimes
\bigwedge^{\bullet}E^*$ 
into the vector space of all diffferential operators
on $\bigwedge^{\bullet} E \otimes \bigwedge^{\bullet} E^*$,
defined on decomposable elements by
\begin{equation}\label{dualembedding}
j_{x \otimes \xi} (y \otimes \eta) = x \wedge i_{\xi} y \otimes \eta \ .
\end{equation}

In an earlier note \cite{B1969}, Buttin had
obtained the Schouten-Nijenhuis bracket of multivectors by a
construction involving an auxiliary torsionless linear connection,~$\nabla$.
For a decomposable element $x$ 
in $V^p(M)$ such that $x= u\wedge v$, $u \in V^{p-1}(M), 
v \in V^1(M)$, set
$\nabla_x y= u \wedge\nabla_v y$, for all $y \in V^{\bullet}(M)$.
To a multivector $x$ of degree $p$, Buttin associated 
the derivation 
of degree $p-1$
of the algebra of multivectors defined by 
\begin{equation}
\widetilde x (y)=  \nabla_x y- j_{\nabla x}y \ , 
\end{equation}
for $y \in
V^{\bullet}(M)$, where $j$ is the map defined by \eqref{dualembedding}.
The Schouten-Nijenhuis bracket of
multivectors $x$ and $y$ is then obtained by letting the derivation
$\widetilde x$ act on $y$.

\subsection{The Vinogradov bracket}\label{vinog}
Vinogradov \cite{V} \cite{CV}
introduced a bilinear operation on the
vector space of all graded endomorphisms of the space of
differential forms on a smooth manifold. (Actually his general
construction is given for any complex.) If $a$ and $b$ are
endomorphisms of the space of differential forms, $\W$, he sets
\begin{equation}\label{eqvinog}
[a,b]_V= \frac{1}{2} \left([[a,{\d}],b] - (-1)^{|b|}[a,[b,{\d}]]\right) \ .
\end{equation}
This bilinear bracket is skew-symmetric but does
not satisfy the Jacobi identity. See \cite{V} for the explicit
trilinear expression of the defect in the Jacobi identity. 
The vector space of all mutivector-valued forms
embeds into this space of
endomorphisms, but it is not closed under this bracket. However, 
the following properties, which will be proved in Section \ref{loday},
are valid.

\medskip

(A) The space of multivectors is closed under the Vinogradov
bracket. Its
restriction to the space of multivectors is a graded Lie bracket, which is 
the Schouten-Nijenhuis bracket. 

\medskip

(B) 
The restriction of the Vinogradov
bracket to the space of vector-valued forms is equal to the
Fr{\"o}licher-Nijenhuis bracket, up to a derivation of
$\Omega^{\bullet}(M)$ of the form $[i_Z,\d]$, where $Z$ is a
vector-valued form.

\medskip

(C) The direct sum of the space of vector fields 
and the space of differential forms is closed under the Vinogradov
bracket. 
The restriction of the 
bracket $[~,~]_V$ to this space was not considered by Vinogradov.
It is skew-symmetric by definition, but it
does not satisfy the Jacobi identity.
\begin{itemize}
\item
Case of $1$-forms: We shall show in Section
\ref{loday} that, when the Vinogradov bracket is 
further restricted to the direct sum of the space of vector fields 
and the space of differential $1$-forms, it is nothing other
than the bracket of Courant \cite{C}.
\item
Case of $p$-forms: 
In fact, the formula for the Courant bracket \eqref{courant} below 
also makes sense in the more general case of a
vector and a form of arbitrary degree. This was observed by Wade in
\cite{W} and independently by Hitchin in
\cite{H}. The calculation in Section \ref{loday}
shows that the bracket defined by
this formula is the restriction of the Vinogradov bracket.
\end{itemize}

\subsection{Unification theorems}
In Buttin and in Vinogradov, we find two ``unification theorems'', 
namely constructions of which both the Schouten-Nijenhuis bracket 
and the Fr{\"o}licher-Nijenhuis bracket
are, in some sense, particular cases. These constructions are
described in different settings: Buttin
introduces a skew-symmetric bilinear 
map from pairs of multivector-valued forms
to differential operators on the space, $\Omega^{\bullet}(M)$, 
of differential forms,
the generalized differential concomitant of the
  first kind,
and she shows that the image of a pair corresponds
to a multivector-valued form only in the case of a pair of
multivectors
or in the case of a pair of vector-valued forms. For his part, 
Vinogradov constructs a skew-symmetric bracket, which does not satisfy the
graded Jacobi identity, defined on pairs of
differential operators on 
$\Omega^{\bullet}(M)$, with values in the space
of differential operators,
which restricts to the Schouten-Nijenhuis bracket 
on multivectors,
and also has the 
property of restricting to the Fr{\"o}licher-Nijenhuis bracket 
on vector-valued forms, but only modulo generalized Lie derivatives.

Using the graded Jacobi identity for the graded commutator, we find
the following relation: 
for any endomorphisms, $a$ and $b$, of $\W$,
\begin{equation}
[[a,b]_V, {\d}] = [[a,{\d}], [b,{\d}]] \ . 
\end{equation}
Therefore, whenever $a = i_X, b=i_{Y}$, for $X$ and $Y$
multivector-valued forms, by \eqref{concom},
\begin{equation}
[\{a,b\}_B,{\d}] = [[a,b]_V, {\d}] \ ,
\end{equation}
and, in particular, if both $X$ and $Y$ are multivectors, or if
both are vector-valued forms,
\begin{equation}
[i_{[X,Y]_B}, {\d}] = [[i_X,i_{Y}]_V,{\d}] \ .
\end{equation}

We claim that the situation can be clarified by the consideration of
non skew-symmetric brackets.

\subsection{Loday brackets on forms and multivectors}\label{loday}

We shall show that the Vinogradov bracket is the skew-symmetrization of
a Loday bracket, which is a derived bracket of the graded commutator of
graded endomorphisms.
In case (A) of Section \ref{vinog}, 
the derived bracket is skew-symmetric and
therefore the Vinogradov bracket coincides with it,
while in cases (B) and (C), 
the derived bracket is not skew-symmetric.
The skew-symmetrization then yields the Vinogradov bracket.

\medskip
 
Again let ${\rm End} (\Omega^{\bullet}(M))$ 
be the algebra of graded endomorphisms of
$\Omega^{\bullet}(M)$, and let $[~,~]$ be the graded commutator.
Let ${\d}$ be the de Rham differential. The derived bracket of
$a, b \in {\rm End} (\Omega^{\bullet}(M))$ is defined by formula
\eqref{interior} of Section \ref{derived},
\begin{equation}
[a,b]_{\d} = [[a,{\d}],b] \ .
\end{equation}
Below we set, for $X \in \Omega ^{\bullet}(M) \otimes V^1(M)$, 
$L_X = [i_X, {\d}]$, where 
$i_X$ is defined by \eqref{embedding}, and is considered as an
endomorphism of $\Omega ^{\bullet}(M) \otimes V^1(M)$.

\begin{theorem}\label{theorem}
(i) The derived bracket $[~,~]_{\d}$
  defines an odd Loday
  algebra structure on ${\rm End}(\W)$.

(ii) The space $V^{\bullet}(M)$ is
closed under the derived bracket $[~,~]_{\d}$. Its restriction to
this subspace of ${\rm End}(\W)$ is skew-symmetric; it is the
Schouten-Nijenhuis bracket.

(iii) The algebraic part of the restricition of $\bra_{\d}$ to 
the subspace $\Omega^{\bullet}(M) \otimes V^{1}(M)$  
of ${\rm End}(\W)$ is the
Fr{\"o}licher-Nijenhuis bracket. More precisely,
\begin{equation}\label{FN}
[i_X,i_Y]_{\d} = i_{[X,Y]_{FN}} - (-1)^{q(q'-1)}L_{i_Y X} \ ,
\end{equation}
where $[X,Y]_{FN}$ denotes the Fr{\"o}licher-Nijenhuis bracket of 
$X \in \Omega ^q(M) \otimes V^1(M)$ and $Y \in \Omega ^{q'}(M) \otimes
V^1(M)$.
  
(iv) The derived brackets of a vector field $x$ and a differerential
form $\xi$ are
\begin{equation}\label{seminonskew}
[x,\xi]_{\d}= L_x\xi 
\end{equation}
and
\begin{equation}\label{seminonskew2}
[\xi,x]_{\d}= -
i_x{\d}\xi \ ,
\end{equation}
and the restriction of  
$[~,~]_{\d}$ to the direct sum of the space of vector fields 
and the space of differential forms is given by
\begin{equation}\label{nonskew}
[x+\xi, y+\eta]_{\d}= [x,y] +  L_x\eta -
i_y{\d}\xi \ ,
\end{equation}
for all vector fields $x$ and $y$, and for all differential forms $\xi$ and
$\eta$, where $[x,y]$ is the Lie bracket of $x$ and $y$.
\end{theorem}  

\noindent{\it Proof} 
Part (i) is a corollary of Theorem \ref{theoremderived} above.
Part (ii) is a re-statement of the Cartan formula \eqref{cartan}
for multivector fields.

To prove part (iii), 
we observe that, for vector-valued forms, $X$ and~$Y$,
$$
[[i_X,i_Y]_{\d},{\d}]=[i_{[X,Y]_{FN}},{\d}] \ ,
$$
because both expressions are equal to $[[i_X,{\d}],[i_Y,{\d}]]$, the
first by the Jacobi identity and the second by definition.
This implies that the algebraic parts of $[i_X,i_Y]_{\d}$ and
$i_{[X,Y]_{FN}}$ are equal.
It follows from formula 5.15 of \cite{FN} 
that the Fr{\"o}licher-Nijenhuis bracket
satisfies,
for $X= \xi \otimes x$, $Y=\eta \otimes y$,
\begin{equation}\label{eq1}
[X,Y]_{FN}= \xi \wedge \eta \otimes [x,y] 
\end{equation}
$$
+ (\xi \wedge L_x \eta + (-1)^{|\xi|} {\d} \xi \wedge i_x\eta) \otimes y
- (-1)^{|\xi||\eta|} 
(\eta \wedge L_y \xi + (-1)^{|\eta|} {\d} \eta \wedge i_y\xi) \otimes x
\ .
$$
Using $i_{\xi \otimes x} =
e_{\xi} \circ i_x$ and $[i_{\xi \otimes x} ,{\d}]= \xi \wedge L_x +
(-1)^{|\xi|} {\d}\xi \wedge i_x$, this formula also implies the expression 
\begin{equation}
[X,Y]_{FN}=  \xi \wedge \eta \otimes [x,y] + L_X\eta \otimes y -
(-1)^{|\xi||\eta|}L_Y \xi \otimes x \ .
\end{equation}
to be found in, {\it e.g.}, \cite{KS} \cite{grabowski2}.

On the other hand, a direct computation shows that
\begin{equation}\label{eq2}
[i_{\xi \otimes x}, i_{\eta \otimes y}]_{\d}=[[i_{\xi \otimes x},{\d}],i_{\eta
  \otimes y}]
\end{equation}
$$
= \xi \wedge \eta \wedge i_{[x,y]} + (\xi \wedge L_x \eta
+ (-1)^{|\xi|} {\d} \xi \wedge i_x \eta)\wedge i _y 
$$
$$
+ (-1)^{|\xi||\eta|+1} \eta \wedge i _y {\d} \xi \wedge i _x 
+ (-1)^{|\xi||\eta|+|\xi|+1} \eta \wedge i _y \xi \wedge L _x \ .
$$
Therefore, we find from \eqref{eq1} and \eqref{eq2},
$$
(-1)^{|\xi||\eta| +|\xi| +1}([i_{\xi \otimes x}, i_{\eta \otimes
  y}]_{\d} 
- i_{[\xi \otimes x, \eta
  \otimes y]_{FN}} )=
 \eta \wedge i_y \xi \wedge L_x - (-1)^{|\xi|+|\eta|} {\d} (\eta \wedge
  i_y\xi) \wedge i_x
$$
$$
= 
 [i_{\eta \wedge i_y\xi \otimes x}, {\d}] =
 [i_{i_YX},{\d}] =
 L_{i_YX}
\ ,
$$
which proves \eqref{FN}. 

To prove part (iv), we first observe that, on $V^1(M)$, 
the derived bracket restricts to
the Lie bracket of vector fields: this follows from the Cartan
formula
\eqref{cartan}. 

If $\xi$ is a differential form of
degree $|\xi|$, 
\begin{equation}
[e_{\xi},{\d}] = (-1)^{|\xi|+1} 
e_{{\d}\xi} \ .
\end{equation} 
Since any two exterior products by forms commute, the derived
bracket vanishes if both arguments are
differential forms.

If $x$ is a vector field and $\xi$ a differential form of
degree $|\xi|$,
\begin{equation}
[i_x,e_{\xi}]_{\d}= [[i_x,{\d}],e_{\xi}]= L_x e_{\xi}-e_{\xi}L_x= 
e_{L_x\xi} \ ,
\end{equation} 
and
\begin{equation}
[e_{\xi},i_x]_{\d}= [[e_{\xi},{\d}],i_x] = (-1)^{|\xi|+1} 
[e_{{\d}\xi},i_x] = - e_{i_x{\d}\xi}.
\end{equation} 
Therefore, identifying vectors and forms with their image under the
embedding $i$, we obtain formulas \eqref{seminonskew} and 
\eqref{seminonskew2}.

To summarize, $\bra_{\d}$
extends the Lie bracket of vector
fields, vanishes on pairs of differential forms and satisfies
\eqref{seminonskew} and \eqref{seminonskew2}, 
and therefore \eqref{nonskew} follows.

\medskip

\noindent{\bf Remark}
The space
$V^{\bullet}(M) \otimes \Omega^{\bullet}(M) \subset 
\End (\Omega^{\bullet}(M))$
is not closed under the derived bracket $\bra_{\d}$ unless the manifold is of
dimension $\leq 1$. To prove this, let us
show that, on manifolds of dimension $\geq 2$, 
there exist multivector fields, $x$ and $y$, 
and differential $1$-forms, $\xi$ and $\eta$, such that  
the operator $[[i_{x \otimes \xi} , {\d}], i_{y\otimes \eta}]$
is not
$C^{\infty}(M)$-linear. For any differential
form 
$\alpha$,
$$
[[i_{\xi \otimes x} , {\d}], i_{\eta \otimes y}](\alpha)
= \xi \wedge L_x {\d}(\eta \wedge i_y\alpha)
- (-1)^{|x|+|\xi|} 
{\d} \xi \wedge i_x(\eta \wedge i_y \alpha)
$$
$$ 
+ (-1)^{|x|+ |y| + |\xi|+|\eta|} \eta \wedge i_y(\xi \wedge L_x
{\d}\alpha)
- (-1)^{|y|+|\eta|} \eta \wedge i_y({\d}\xi \wedge i_x \alpha) \ .
$$
If $x$ and $y$ are bivectors, 
for a $1$-form $\alpha$ and a function $f$, 
$$
 [[i_{\xi \otimes x} , {\d}], i_{\eta \otimes y}](f \alpha)
- f [[i_{\xi \otimes x} , {\d}], i_{\eta \otimes y}] (\alpha)
$$
$$
= \left(i_y(\xi \wedge i_x({\d}f \wedge {\d}\alpha)) -
  (-1)^{|x|}(i_x{\d}\alpha) i_y(\xi\wedge {\d}f)\right)\eta \ .
$$
If ${\d}\alpha =\beta \wedge \gamma$, this expression is 
$$
 \left( i_x(\beta \wedge {\d}f) i_y (\xi \wedge \gamma) - 
i_x (\gamma \wedge {\d}f)
i_y (\xi \wedge\beta) \right) \eta \ ,
$$
which does not vanish in general, as can be proved by using
local coordinates.
 
\begin{theorem}
(i) The skew-symmetrization of the derived bracket
$[~,~]_{\d}$
is the Vinogradov bracket.

(ii) The skew-symmetrized derived bracket of a vector $x$ and a form $\xi$
is
\begin{equation}
[x,\xi]^{-}_{\d}
 = L_x\xi -\frac{1}{2} {\d}i_x \xi \ , 
\end{equation}
and the restriction of the skew-symmetrization of 
$[~,~]_{\d}$ to the direct sum of the space of vector fields 
and the space of differential forms is
given by 
\begin{equation}\label{courant}
[x+\xi,y+\eta]^{-}_{\d} = 
[x,y] + L_x\eta - L_y\xi  -\frac{1}{2} {\d}(i_x \eta -
i_y \xi) \ ,
\end{equation}
for all vector fields $x$ and $y$, and all differential forms $\xi$ and
$\eta$.

(iii) The restriction of the skew-symmetrization of 
$[~,~]_{\d}$ to the direct sum of the space of vector fields 
and the space of differential $1$-forms is
the Courant bracket.
\end{theorem}

\noindent{\it Proof}
Part (i) follows from equations \eqref{skewint} and \eqref{eqvinog}.
Part (ii) follows immediately from part (iv) of Theorem \ref{theorem}.
On $V^1(M) \oplus
\Omega^1(M)$, formula \eqref{courant} is precisely the Courant bracket 
as defined in \cite{C} and used in \cite{LWX}, proving part (iii). 

\medskip

\noindent{\bf Remark}
If we consider the restriction of the skew-symmetrized derived bracket
to the vector-valued forms, we find
$$
[i_X,i_Y]_V=[i_X,i_Y]^-_{\d}= i_{[X,Y]_{FN}} + \frac{1}{2} (-1)^p
L_{i_XY + (-1)^{(p-1)(p'-1)}i_YX} \ .
$$
It is clear that $\Omega^{\bullet}(M) \otimes V^1(M)$
is not closed under the derived bracket nor under its
skew-symmetrization, because neither 
$[i_{\xi \otimes x},i_{\eta \otimes y}]_{\d}$ nor
$[i_{\xi \otimes x},i_{\eta \otimes y}]^-_{\d}$ vanishes on
functions. In fact,
$$
[i_{\xi \otimes x},i_{\eta \otimes y}]_{\d}f = \eta \wedge i_y\xi L_xf
\ .
$$

\medskip

The explicit expression \eqref{nonskew} of the derived bracket 
for the case of $1$-forms appears in
Dorfman (see \cite{Dorf87}, 
\cite{Dorf}), 
in her study
of the properties of Dirac structures on complexes over Lie algebras.
However, that expression 
does not appear as a derived bracket and its properties
are not spelled out.

The relationship of the derived bracket with the Vinogradov bracket
was shown in \cite{yks96}, but its relationship with the Courant
bracket was observed later,
independently by myself (in an e.mail letter to Alan Weinstein, 1998), 
Pavol \v Severa and Ping Xu 
(all unpublished, see \cite{SW}). 
It is the non
skew-symmetric bracket which 
is now used in the theory of Courant algebroids
\cite{Rphd} \cite{Rlmp} 
\cite{SW}.

\section{Odd and even brackets on supermanifolds}
In the interpretation of the brackets of differential geometry in
terms of supermanifolds, the odd brackets are obtained as derived
brackets of even brackets and conversely. 
This was shown in all generality by Voronov in
\cite{Vo}. 
See also the article by Batalin and Marnelius \cite{BM}.
The main examples are 

(1) the Schouten-Nijenhuis bracket of multivectors on a manifold,
as a derived bracket of the
canonical Poisson bracket on the cotangent bundle 
of the manifold,

(2) the Poisson bracket of functions on a manifold, $M$, as a derived
    bracket of the canonical Schouten-Nijenhuis bracket on the
    cotangent bundle of $M$ with reversed parity (see \cite{yks96}),

(3) the algebraic Schouten bracket on the exterior algebra of a Lie
algebra $(F, \mu)$,
as the derived bracket of the canonical Poisson structure on
$\bigwedge(F\oplus F^*)$ (see \cite{yksJ}).
 
All these instances are particular cases of the general construction 
on Lie algebroids which we shall describe in Section \ref{algd}.

Below, a function on $T^*\M$, where $\M$ is a supermanifold, 
is called a {\it hamiltonian} on $\M$,
and 
a bracket on the vector space $C^{\infty}(\M)$ is sometimes
called a bracket on $\M$.
If $E \to \M$ is a vector bundle, $\Pi E$ denotes the supermanifold
obtained by reversing the parity of the fibers.

\subsection{Odd and even 
Poisson brackets on supermanifolds are derived brackets}
The following theorems are due to
Voronov \cite{Vo} 
(see also \cite{BM}). 

\begin{theorem}\label{poisson}
 Any odd Poisson bracket on 
a supermanifold, $\M$,  
is a derived
  bracket of the canonical Poisson bracket $\{~,~\}$ on $T^*\M$. 
More precisely, for any odd Poisson bracket on $\M$,
$[~,~]$, there exists a quadratic function $S$ on $T^*\M$ such
that 
\begin{equation}\label{hamilt}
[f,g] = \{\{f,S\},g\} \ ,
\end{equation}
for all $f, g \in C^{\infty}(\M)$.
(On the right-hand side of \eqref{hamilt}, $f$ and $g$
are identified with their pull-backs to $T^*\M$, {\it i.e.}, are
considered as functions on the
vector bundle
$T^*\M$,
that are constant on the fibers.)
\end{theorem}

Formula \eqref{hamilt} defines a derived bracket $\bra =\{~,~\}_S$ 
which is not only a Loday bracket but is a true Lie bracket, 
{\it i.e.}, it is  
skew-symmetric, because $C^{\infty}(\M)$ is an abelian subalgebra of the
Poisson algebra of $T^*\M$.

In coordinates, 
if $(x^\alpha)$ are local coordinates on $\M$ and $(x^{\alpha}, p_{\alpha})$
the asssociated 
local cordinates on $T^*\M$, then $S=\frac{1}{2} S^{\alpha
  \beta}(x)p_{\alpha} p_{\beta}$, where
$S^{\alpha \beta} = [x^{\alpha},x^{\beta}]$.

\medskip

\noindent{\bf Example} As an exercise, 
let us illustrate this theorem in the
case where $\M =
\Pi T^*N$, for $N$ a manifold of dimension $n$.
The Schouten-Nijenhuis bracket of fields of multivectors is a
canonically defined odd Poisson bracket, $[~,~]_{SN}$, on $\M$.  
What is the corresponding quadratic hamiltonian?
We seek its expression in local coordinates. Let$(y^i,\tilde y_i)$
be adapted local coordinates on $\M = \Pi T^*N$ ($y^i$ is even and  
$\t y_i$
is odd). Then $[y^i,y^j]_{SN}=0$, $[y^i,\t y_j]_{SN}=\delta^i_j$ and 
$[\t y_i,\t y_j]_{SN}=0$. Let $(y^i,\t y_i,p_i,\t p^i)$ 
be associated local coordinates on $T^*\M$. The canonical Poisson bracket
satisfies,
$\{y^i,p_j\}=\delta^i_j$, $\{\t y_i,\t p^j\}=\delta^j_i$, while all
other brackets vanish.
The quadratic hamiltonian $S = - p_i \t p^i\in C^{\infty}(T^*(\Pi T^*N))$ 
satisfies equation \eqref{hamilt} above in the form
$$
[f,g]_{SN}= \{\{f,S\},g\} \ .
$$
The hamiltonian $S$ can be also defined in an invariant way.
For each $\varphi \in T^*\M$, 
$S(\varphi) = - < p^*_{\Pi T^*N}\varphi, p^*_{\Pi TN} (\kappa \varphi)>$,
where $ p^*_{\Pi T^*N}$ 
maps $T^*(\Pi T^*N)$ to $\Pi T^*N$, and $p^*_{\Pi TN}$ maps
  $T^*(\Pi TN)$ to $\Pi TN$, while $\kappa$ is the canonical isomorphism
\cite{Vo} from $T^*(\Pi T^*N)$ to $T^*(\Pi TN)$. 

\bigskip

In the same manner, we can describe every even Poisson bracket on $\M$
as a derived bracket.
Let $P$ be the even Poisson
bivector defining an even Poisson bracket $\{~,~\}$ on $\M$.
Then 
\begin{equation*} 
\{f,g\} = [[f,P],g] \ ,
\end{equation*}
where the bracket on the right-hand is the Schouten-Nijenhuis bracket
of multivector fields on $\M$, and $f$ and $g$ are considered as 
multivector fields of
degree $0$.
See formula (3.5) of \cite{yks96} for the case of ordinary
manifolds, and \cite {Vo} for supermanifolds. 
Since a multivector field on $\M$ is a function on $\Pi T^*\M$, 
this property can be reformulated in the language of supermanifolds,
making it ``dual'' to Theorem \ref{hamilt}.

\begin{theorem}\label{schouten}
 Any even Poisson bracket on $\M$ is a derived
  bracket of the canonical Schouten-Nijenhuis
 bracket $[~,~]$ on $\Pi T^*\M$. 
More precisely, for any even Poisson bracket on $\M$,
$\{~,~\}$, there exists a quadratic function $P$ on $\Pi T^*\M$
such that 
\begin{equation} \label{biv}
\{f,g\} = [[f,P],g] \ ,
\end{equation}
for all $f, g \in C^{\infty}(\M)$. (On the right-hand side, $f$ and $g$
are identified with their pull-backs to $\Pi T^*\M$, {\it i.e.}, 
are considered as functions on the vector 
bundle $\Pi T^*\M$, that are constant on the fibers.)
\end{theorem}

The quadratic function $P$ on $\Pi T^*\M$ is nothing but the Poisson
bivector giving rise to the given even Poisson bracket. 
In local coordinates $(x^{\alpha}, \t p_{\alpha})$ on $\Pi T^*\M$, 
$P=\frac{1}{2} P^{\alpha \beta}(x) \t p_{\beta} \t p_{\alpha}$,
where $ P^{\alpha
  \beta} = \{x^{\alpha}, x^{\beta}\}$.

\subsection{Derived brackets and Lie algebras}\label{bigb}

Let $E$ be a finite-dimensional vector space over a field of
characteristic $0$. For simplicity, we consider the ungraded (purely
even) case, but
the properties below extend to the case where the vector space $E$
is itself graded, see \cite{Vo}.
The following structures are equivalent:
\begin{itemize}
\item{
a {\it Lie algebra structure} on $E$,}
\item{
a {\it Lie co-algebra structure} on $E^*$,}
\item{
a {\it linear Schouten structure} on $\Pi E^*$, {\it i.e.},
a Schouten algebra structure
on $C^{\infty}(\Pi E^*) 
= \bigwedge^{\bullet} E$ such that $E$
is closed under the bracket,}
\item{
a {\it linear Poisson structure} on $E^*$, {\it i.e.}, 
a Poisson algebra structure
on $C^{\infty}(E^*)$ such that $E$
is closed under the bracket, {\it i.e.}, a linear bivector field on $E^*$,}
\item{
a {\it linear-quadratic hamiltonian of Poisson square} $0$ on $\Pi
E^*$, {\it i.e.},
an element $\mu \in C^{\infty}(T^*(\Pi E^*)) =  C^{\infty}(\Pi E^*
\oplus \Pi E)
= \bigwedge^{\bullet}(E\oplus E^*)$ such that $\mu \in \bigwedge^2E^*
\otimes E$ and $\{\mu, \mu\} =0$, also denoted by $H$,}
\item{
a {\it quadratic homological vector field} on $\Pi E$, {\it i.e.},
a quadratic differential on $C^{\infty}(\Pi E) =
\bigwedge^{\bullet}E^*$, often denoted by $d$, or $d_{\mu}$, or $Q$.}
\end{itemize}

The canonical Poisson bracket on 
$C^{\infty}(T^*(\Pi E^*))
= \bigwedge^{\bullet}(E\oplus E^*)$ 
was first defined in \cite{KosSte} and considered in \cite{LR}. It was
an esential tool in \cite{yksJ}, where we first 
called it the {\it big bracket}.
Here, we have denoted the big bracket 
by $\{~,~\}$. 

The bracket on $C^{\infty}(\Pi E^*)= \bigwedge^{\bullet} E$ is
called the {\it algebraic Schouten bracket} of the Lie algebra.
The differential on $C^{\infty}(\Pi E)= \bigwedge^{\bullet}E^*$, 
corresponding to the
Lie algebra structure on $E$, 
is the {\it Chevalley-Eilenberg differential} on the scalar-valued
cochains on $E$.

Let $(e_i)$ be a basis of $E$, with coordinates $(x^i)$, and set 
$\mu (e_i,e_j)=C^k_{ij}e_k$.
Let $(\t \xi_i)$ be the coordinates in the dual basis on $\Pi E^*$, and let 
$(\t \xi_i,
\t x^i)$ be the associated coordinates on $T^*(\Pi E^*)= \Pi E^* \oplus
\Pi E$.
The hamiltonian on $\Pi E^*$, which is a function on
$T^*(\Pi E^*)$, may be written
$$
H =
\frac{1}{2} C_{ij}^k \t \xi_k \t x^i \t x^j \ , 
$$
while the vector field on $\Pi E$ may be written 
$$
Q =
\frac{1}{2} 
\t x^j \t x^i C^k_{ij} \frac{\partial}{\partial \t x^k} \ .
$$

Since Lie algebra structures on $E$ are Poisson (resp., Schouten)
structures on $E^*$ (resp., $\Pi E^*$), 
Theorems \ref{poisson} and \ref{schouten} apply. They
take the following form in the case of Lie algebras.

\begin{corollary}
Given a Lie algebra structure $\mu
\in \bigwedge^2E^*
\otimes E$ on a vector space, $E$, 

(i) the Schouten bracket on 
$C^{\infty}(\Pi E^*)= \bigwedge^{\bullet} E$ is given by the derived bracket
formula,
\begin{equation}\label{alghamilt}
[x,y]_{\mu}= \{\{x,\mu\},y\} \ ,
\end{equation}
where $\{~,~\}$ denotes the canonical Poisson bracket (big bracket)
on 
$T^*(\Pi E^*)$, $\mu$ is considered as a hamiltonian on $\Pi E^*$,
{\it i.e.}, as a function on 
$T^*(\Pi E^*)$, 
and $x$ and $y$ are considered as functions on
$T^*(\Pi E^*)$ that are constant on $\Pi E$,

(ii) the Poisson bracket of $f$ and $g$ in $C^{\infty} (E^*)$ is given by the
derived bracket formula,
\begin{equation}\label{algbiv}
\{f,g\}_{\mu} = [[f, \mu],g] \ ,
\end{equation}
where $[~,~]$ denotes the canonical 
Schouten-Nijenhuis bracket of multivector fields on
$E^*$, 
$\mu$ is considered as a bivector field on $E^*$,
and $f$ and $g$ are considered as multivector fields of degree $0$ on $E^*$.
\end{corollary}

Part (i) of the corollary was observed
in \cite{R} as well as in \cite{yksJ}, where it was used to prove
various properties of Lie bialgebras ands Poisson Lie groups.
As an application, we now recall how to derive 
the condition for an $r$-matrix to define
a coboundary Lie  bialgebra.

\medskip

\noindent{{\bf Example}} Let $\mathfrak g = (E, \mu)$ be a Lie
algebra. If $r \in \bigwedge^2 \mathfrak g$,
by \eqref{alghamilt}, the algebraic Schouten bracket, $[r,r]_{\mu}$, satisfies
$[r,r]_{\mu}=\{\{r,\mu\},r\}$, where $\{~,~\}$ denotes the big bracket. 
Let $d_{\mu}r$ be the Chevalley-Eilenberg
coboundary of $r$. In order for $d_{\mu}r$ to be a Lie cobracket on $E$
it is necessary and sufficient that $\{d_{\mu}r,d_{\mu}r\}=0$.
Using the relations $d_{\mu}r = \{\mu,r\}$ and $\{\mu,\mu\}=0$,
the Jacobi identity and
\eqref{alghamilt}, we obtain
$$
\{d_{\mu}r,d_{\mu}r\}= \{\{\mu,r\},\{\mu,r\}\} = 
\{\mu,\{\{r,\mu\},r\}\}= \{\mu, [r,r]_{\mu}\}= d_{\mu}[r,r]_{\mu} \ .
$$
Therefore 
$d_{\mu}r$ is a Lie cobracket on $E$
if and only if $[r,r]_{\mu}$ is ad-invariant.
Since the Drinfeld bracket,
$<r,r>$, coincides with the algebraic Schouten bracket up to 
a factor $- \frac{1}{2}$, the above computation is a short proof of
the fact that $(E,\mu,d_{\mu}r)$ is a
coboundary Lie
bialgebra if and only if $r$ satisfies 
the {\it generalized classical Yang-Baxter
equation}, {\it i.e.}, the ad-invariance of the 
Drinfeld bracket $<r,r>$.

\medskip

We now state another derived bracket formula in the theory of Lie
algebras
(see \cite{Vo}).
We can consider $x \in E$ as a constant vector field on $\Pi E$.
Let $i: x \in  E
\mapsto i_x \in V^{1}(\Pi E)= {\rm Der}(C^{\infty}(\Pi E)) = {\rm Der}
(\bigwedge^{\bullet}E^*)$
be the canonical embedding. With the preceding notations,
for $x$ and $y$ in 
$E$,
\begin{equation}\label{liealg}
i_{[x,y]_{\mu}} =[[i_x,d_{\mu}],i_y] \ ,
\end{equation}
where the bracket on the right-hand side is the graded commutator.
If, for example, $\alpha$ is a $1$-form on $E$,
this formula reduces to
$$
(d_{\mu} \alpha)(x,y) = - \alpha([x,y]_{\mu}) \ .
$$
More generally 
for $x \in \bigwedge^{\bullet}E$, let $i_x \in
{\rm End}(\bigwedge^{\bullet}E^*)$
be the interior product by $x$. Formula \eqref{liealg} is then valid 
for $x$ and $y \in \bigwedge^{\bullet}E$,
where the bracket on the left-hand side is the algebraic Schouten
bracket.

The axioms of Lie bialgebras and generalizations thereof 
can be easily formulated in this
framework, as shown in \cite{LR}
\cite{yksJ} \cite{Rphd} 
\cite{R2002} \cite{Vo}.

\medskip

In conclusion, we can state: just as the
Lie bracket of vector fields is a derived bracket according to
the Cartan relation \eqref{cartan}, the Lie bracket on any Lie algebra is a
derived bracket according to equation \eqref{liealg}.
We observe that the de Rham differential and the Chevalley-Eilenberg
cohomology operator play analogous roles.
The Lie algebroid framework which we shall now describe unifies 
these two theories.

\subsection{Derived brackets and Lie algebroids}\label{algd}
The approach to Lie algebroids in terms of supermanifolds is due to
Vaintrob \cite{Vaintrob}, and was developped by Roytenberg \cite{Rphd}
\cite{Rlmp} \cite{R2002} and by Voronov \cite{Vo}. 
See also
Alexandrov, Schwarz, Zaboronsky and Kontsevich 
\cite{ASZK} and Jae-Suk Park \cite{P}.

Let $A \to M$ be a vector bundle.
A {\it Lie algebroid structure} on $A$ can
be defined in several equivalent ways:
\begin{itemize}
\item{
a {\it Lie algebroid structure} on $A$, {\it i.e.}, a Lie algebra structure
on $\Gamma A$ and a morphism of vector bundles, $\rho : A \to TM$,
called the anchor,
satisfying the Leibniz rule, 
\begin{equation}
[u,fv]=f[u,v] +(\rho(u)f) v \ .
\end{equation}
for all $u$ and $v \in \Gamma A$, $f \in C^{\infty}(M)$,}
\item{
a {\it linear Schouten structure} on $\Pi A^*$, {\it i.e.}, 
a Schouten algebra structure
on $C^{\infty}(\Pi A^*)= \Gamma(\bigwedge^{\bullet} A)$ such that
$\Gamma A$
is closed under the bracket,}
\item{
a {\it linear Poisson structure} on $A^*$, {\it i.e.}, 
a Poisson algebra structure
on $C^{\infty}(A^*)$ such that $\Gamma A$
is closed under the bracket.}
\end{itemize}

\noindent{\bf Remark}
The fact that the anchor maps the bracket on $\Gamma A$
to the Lie bracket
on $\Gamma TM$, which is usually listed in the axioms of a Lie algebroid,
is actually a consequence of the Jacobi identity \eqref{jac}
for the bracket on
$\Gamma A$ together with the Leibniz rule, as we show by computing
$[u,[v,fw]]$ in two ways (see \cite{yksm}).
\medskip

A Lie algebroid structure on $A$ 
can also be defined by a {\it homological vector field} 
on $\Pi A$, {\it i.e.},
a differential on $C^{\infty}(\Pi A) =
\Gamma(\bigwedge^{\bullet}A^*)$, the {\it Lie algebroid differential},
often denoted by $Q$, or $d$, or $d_A$.
Let $(x^{\alpha})$ be local coordinates on $M$, let
$(e_{i})$ be a local basis of $\Gamma A$, 
and let $(x^{\alpha}, y^i)$ be the corresponding local coordinates on $A$.
Let
$$
\rho(e_i) = a_i^{\alpha}(x) 
\frac{\partial}{\partial x^{\alpha}}
$$
and 
$$ 
[e_i,e_j] = C^k_{ij}(x) e_k \ .
$$
Then the vector field, $Q$,  on the vector bundle $\Pi A$, equipped 
with the local coordinates $(x^{\alpha}, \t y^i)$, has the local expression,
$$
Q = \t y^i a_i^{\alpha}(x) \frac{\partial}{\partial x^{\alpha}}
+ \frac{1}{2} \t y^j \t y^i C^k_{ij}(x) \frac{\partial}{\partial
  \t y^{k}} \ .
$$

A Lie algebroid structure on $A$ can also be viewed as 
a {\it quadratic hamiltonian} on $\Pi A^*$, of Poisson square $0$, denoted
by $H$. 
If $(x^{\alpha},\t \eta_i)$ 
are the local coordinates on
$\Pi A^*$ dual to $(x^{\alpha},\t y^i)$, 
and $(x^{\alpha},\t \eta_i, p_{\alpha}, \t \theta^i)$ are the
associated local coordinates on $T^*(\Pi A^*)$, then
$$
H = 
a_i^{\alpha}(x) p_{\alpha} \t \theta^i 
+  \frac{1}{2}  \t \eta_k
C^k_{ij}(x) \t \theta^j \t \theta^i \ .
$$

The structure can also be viewed as a {\it bivector field} on 
$A^*$, defining the linear Poisson structure, denoted by $P$. 
Its local expression in
the local coordinates $(x^{\alpha}, \eta_i)$ on $A^*$ is 
$$
P =  a_i^{\alpha}(x) \frac{\partial}{\partial x^{\alpha}}
\frac{\partial}{\partial \eta_i} 
+  \frac{1}{2}  \eta_k
C^k_{ij}(x) \frac{\partial}{\partial \eta_j}
\frac{\partial}{\partial  \eta_i} \ .
$$

There are several derived brackets in this theory.
\begin{theorem}
Given a Lie algebroid structure on the vector bundle~$A$,

(i) if $H$ is the hamiltonian of the Lie algebroid, then the Schouten
bracket of $u$ and $v \in C^{\infty}(\Pi A^*) = 
\Gamma(\bigwedge^{\bullet}A)$ is the derived bracket,
\begin{equation}\label{LAhamilt}
[u,v]_A= \{\{u,H\},v\} \ ,
\end{equation}
where $\{~,~\}$ is the canonical Poisson bracket of $T^*(\Pi A^*)$,
and $u$ and $v$ are considered as functions on $T^*(\Pi A^*)$
that are constant on the fibers,

(ii) if $P$ is the bivector field 
on $A^*$ of the Lie algebroid, then
    the Poisson bracket of $\varphi$ and $\psi \in C^{\infty}(A^*)$ is
the derived bracket,
\begin{equation}\label{LAbiv}
\{\varphi, \psi\}_A= [[\varphi, P], \psi],
\end{equation}
where $[~,~]$ is the canonical Schouten-Nijenhuis 
bracket of multivector fields on 
$A^*$, and $\varphi$
and $\psi$ are considered as multivectors of degree $0$.
\end{theorem}

In particular, formula \eqref{LAhamilt} is valid for $u$ and $v$ sections of
$A$, showing that any Lie algebroid bracket is a derived bracket
of a canonical Poisson bracket by a quadratic hamiltonian.

\medskip

In terms of endomorphisms of $C^{\infty}(\Pi A)=
\Gamma(\bigwedge^{\bullet}A^*)$,
\begin{equation}\label{LAend}
i_{[u,v]_A}=  [[i_u,d_A],i_v] \ ,
\end{equation}
where $u$ and $v$ are functions on $\Pi A^*$, $d_A$ is the Lie
algebroid differential, and $i$ is the interior
product.
When $u$
and $v$ are sections of $A$, relation \eqref{LAend} 
is an equality of derivations of  
$\Gamma(\bigwedge^{\bullet}A^*)$.

As a particular case of \eqref{LAend}, we see once more that the Lie
algebroid bracket is a derived bracket. Moreover, if $u \in \Gamma A$
and $f \in C^{\infty}(M)$,
$[u,f]_A$ is the function $[[i_u,d_A],f]= \rho(u)f$.
Thus, the action of the anchor also appears as a derived bracket. 
\medskip

Also, on a Lie algebroid, $A$, a
Fr{\"o}licher-Nijenhuis bracket, $\bra_{FN}$, can be defined \cite{DVM}
\cite{grabowski2}
on vector-valued forms, {\it i.e.},
sections of $\bigwedge^{\bullet}A^* \otimes A$, in such a way that 
it satisfies
\begin{equation}
[i_{[X,Y]_{FN}}, d_A] =[[i_X,d_A],
[i_Y, d_A]] \ .
\end{equation}
In fact, this generalized Fr\"olicher-Nijenhuis bracket is defined by
means of equation \eqref{eq1},
in which the de Rham differential and the Lie derivation are replaced by
their Lie algebroid generalizations.

\medskip

Obviously, the formulas for brackets on Lie algebroids admit two
particular cases of interest:
\begin{itemize}
\item Lie algebras. When the base manifold is a point, we recover the
  case of a Lie algebra. Formula \eqref{LAhamilt} 
reduces to formula \eqref{alghamilt}, and formula \eqref{LAbiv} 
reduces to formula \eqref{algbiv}, 
while formula \eqref{LAend} reduces to \eqref{liealg}.
\item Manifolds. When the Lie algebroid $A$ is a tangent bundle,
  we recover the case of a manifold. The differential $d_A$ is the de
  Rham differential.
Formula \eqref{LAhamilt} reduces to formula 
\eqref{hamilt}, and formula \eqref{LAbiv} reduces to formula \eqref{biv},
while formula \eqref{LAend} reduces to the Cartan formula
\eqref{cartan}.
\end{itemize}

As an application of the general formulas, 
we mention the case $A=T^*M$, when $M$ is a
Poisson manifold with Poisson bivector $P$. Then it is known (see 
\cite{yksm}) that
$d_A$ is the Lichnerowicz-Poisson differential, 
$d_P=[P, ~ \cdot ~ ]_{SN}$, where
$[~,~]_{SN}$ is the Schouten-Nijenhuis bracket of multivectors.
If $\bra^P$ denotes the Koszul bracket of differential forms,
then, by formula \eqref{LAend},
\begin{equation}\label{kras}
i_{[\alpha, \beta]^P}=[[i_{\alpha},d_P ],i_{\beta}] \ ,
\end{equation}
for any forms $\alpha$ and $\beta$. This formula
was first proved by
Krasilsh'chilk
in \cite{Kr}. (See also \cite{yks96}.) It appears once more that
Poisson structures play a role dual to that of differential structures
on manifolds: in this sense, formula \eqref{kras} is dual to the
Cartan formula \eqref{cartan}.

\section{Derived brackets and Courant algebroids}\label{courantalg}

Courant algebroids were introduced by Liu, Weinstein and Xu in
\cite{LWX}, as a generalization of the bracket
defined by Courant \cite{C} on the sections of $TM \oplus T^*M$. As we
pointed out in Sections \ref{vinog} and \ref{loday}, 
the Courant bracket is
skew-symmetric but
does not satisfy the Jacobi identity. The definition of Courant
algebroids
was later
re-formulated \cite{Rphd} \cite{SW} \cite{R2002} in terms of Loday brackets.
An in-depth study of the role of derived brackets
in the general theory of Courant algebroids can be found in the article by
Roytenberg
\cite{R2002}, in terms of cubic hamiltonians on graded
supermanifolds, and in 
the forthcoming study by Alekseev and Xu \cite{AX}, in terms of
Clifford modules and compatible connections.
Here we shall be content with exhibiting one more instance of a
derived bracket construction.

\subsection{Courant algebroids}
The vector bundle $TM \oplus T^*M$, for any manifold $M$, 
with the field of nondegenerate symmetric bilinear forms
defined by the conditions that $TM$ and $T^*M$
be isotropic and $(x|\xi)=<x,\xi>$, for $x$ a tangent vector and $\xi$
a $1$-form at a point in $M$,
with the bracket
\eqref{nonskew} and anchor the projection onto $TM$, is the
prototypical exemple of a Courant
algebroid.
As shown in \cite{LWX}, this construction 
can be generalized to the case of $A \oplus A^*$, where
$(A,A^*)$ is a {\it Lie bialgebroid}, {\it i.e.}, 
a pair of Lie algebroids in duality
satisfying a compatibility condition \cite{MX} \cite{yksbialg}
\cite{Rlmp}: taking into account the bracket of $A^*$, 
formula \eqref{nonskew} can be extended to the sections of $A
\oplus A^*$, and together with the sum of the anchors of $A$ and
$A^*$, it defines a Courant algebroid structure on $A \oplus A^*$.

\subsection{Courant bracket with background as a derived bracket}
As an example, let us find an explicit deriving operator for the
Courant algebroid associated to a
Poisson structure with background in the sense of
\v Severa and Weinstein 
\cite{SW}, {\it i.e.}, 
a WZW-Poisson structure in the sense of Klim\v cik and
Strobl \cite{KSt}.

Let $\psi$ be an
arbitrary form of odd degree on a manifold $M$.
We consider the operator on $\Omega^{\bullet}(M)$,
$$
{\d}^{\psi} = {\d} +e_{\psi} \ .
$$
Then $[{\d}^{\psi},{\d}^{\psi}]=
e_{{\d}\psi}$.
So, whenever $\psi$ is a {\it closed} form of odd degree, we can consider
the derived bracket on $\End(\Omega^{\bullet}(M))$, 
$\bra_{{\d}^{\psi}}$, arising from the graded
commutator and the odd interior derivation of square $0$ defined by
${\d}^{\psi}$. For vector fields $x$ and $y$,
$$
[i_x,i_y]_{{\d}^{\psi}} = [[i_x,{\d}^{\psi}],i_y] =
i_{[x,y]} + e_{i_{x\wedge y}\psi}, 
$$
where $[x,y]$ is the Lie bracket of $x$ and $y$.
The first non-trivial case is when $\psi$ is of degree $3$.
We see that $V^1(M) \oplus \Omega^1(M)$ is closed under the
derived bracket $\bra_{{\d}^{\psi}}$ if and only if $\psi$ is a {\it form of
degree}~$3$. 
Therefore, let $\psi$ be a {\it closed} $3$-{\it form}.
We define
$$
[x+\xi, y +\eta]_{{\d}^{\psi}}= [[x+\xi,{{\d}^{\psi}}],y+{\eta}] \ .
$$
for vector fields $x$ and $y$, and $1$-forms $\xi$ and
$\eta$,
and we find the following generalization of \eqref{nonskew},
\begin{equation}
[x+\xi, y +\eta]_{{\d}^{\psi}}= [x,y] + L_x\eta -i_y {\d}\xi +i_{x
  \wedge y}\psi \ .
\end{equation}
In particular, the bracket of any two differential $1$-forms remains $0$,
but the bracket of two vector fields has both a component in the
space of vector fields and a component in the space of differential
$1$-forms.

This bracket, together with the field of symmetric bilinear forms
recalled above, and, for anchor, the projection onto $TM$,
turn $TM \oplus T^*M$ into a Courant algebroid, called the 
{\it Courant algebroid with background} $\psi$ \cite{SW}. 
We have just shown that it
is a derived bracket by a modified differential.
 
\begin{proposition}
The 
Courant bracket with background $\psi$
on $TM \oplus T^*M$ is the derived bracket of the commutator of
endomorphisms of
$\Omega ^{\bullet}(M)$ by ${\d} + e_{\psi}$.
\end{proposition}  

\subsection{Properties of Courant brackets with background}
To conclude, we list a few properties of the Courant
algebroids with background, following mainly \cite{SW}. 
Now let $P$ be a bivector on $M$, and let 
$P^{\sharp}$ be the mapping from $T^*M$
to $TM$, defined by $P^{\sharp}\xi = i_{\xi}P$.
We shall determine the condition for the graph of $P^{\sharp}$ to be 
a $\psi$-{\it Dirac structure}, {\it i.e.}, to be maximally isotropic
and closed 
under the bracket $\bra_{{\d}^{\psi}}$.
Let $\xi$ and $\eta$ be $1$-forms.
Then
$$
[P^{\sharp}\xi + \xi, P^{\sharp}\eta + \eta ]_{{\d}^{\psi}}
= [P^{\sharp}\xi,P^{\sharp}\eta] + L_{P^{\sharp}\xi}\eta -
i_{P^{\sharp}\eta}
{\d}\xi+ i_{P^{\sharp}\xi \wedge P^{\sharp}\eta}\psi 
$$
$$
= [P^{\sharp}\xi,P^{\sharp}\eta] + L_{P^{\sharp}\xi}\eta -
L_{P^{\sharp}\eta}\xi - {\d}(P(\xi,\eta)) 
+ i_{P^{\sharp}\xi \wedge P^{\sharp}\eta}\psi 
$$
$$
= [P^{\sharp}\xi,P^{\sharp}\eta] + [\xi,\eta]^P 
+ i_{P^{\sharp}\xi \wedge P^{\sharp}\eta}\psi \ , 
$$
and therefore the condition for the graph of $P^{\sharp}$ 
to be closed under the
derived bracket is that
$$
P^{\sharp}([\xi,\eta]^P 
+ i_{P^{\sharp}\xi \wedge P^{\sharp}\eta}\psi) -  [P^{\sharp}\xi ,
P^{\sharp}\eta] = 0 \ .
$$
This condition is equivalent to 
\begin{equation} \label{cns}
\frac{1}{2} [P,P]_{SN} = (\wedge^3 P^{\sharp}) (\psi)  \ ,
\end{equation}
where $[P,P]_{SN}$ is the Schouten-Nijenhuis bracket of $P$ with itself.
If the graph of $P^{\sharp}$ is a
$\psi$-Dirac structure, {\it i.e.}, if condition \eqref{cns} is satisfied, 
$P$ is called a {\it Poisson
structure with background} $\psi$. 

We set
\begin{equation}\label{bracket}
[\xi,\eta]^{P,\psi} = [\xi,\eta]^P
+i_{P^{\sharp}\xi \wedge P^{\sharp}\eta}\psi  \ .
\end{equation}
Formula \eqref{bracket}
defines a skew-symmetric bracket on $T^*M$ with anchor $P^{\sharp}$.
(The Lie bracket of $1$-forms
defined by a Poisson bivector is recovered as the special case
$\psi = 0$.)

The corresponding derivation 
$d_{P,\psi}$ on  $V^{\bullet}(M)$ 
satisfies
$d_{P,\psi} f = 
d_{P}f$,
for $f \in \C$, and 
$$
(d_{P,\psi} x)(\xi, \eta) = 
P^{\sharp}\xi <\eta,x> - \,
P^{\sharp}\eta <\xi, x> - < [\xi,\eta]^{P,\psi},x> \ ,
$$
for all $x \in V^1(M)$, $\xi, \eta \in \Omega^1(M)$.
Let us 
define $(\wedge^2 P^{\sharp}) (\psi)$ to be the bivector-valued $1$-form
such that $(\wedge^2 P^{\sharp})(\psi)(x)(\xi,\eta) 
= \psi(P^{\sharp}\xi, P^{\sharp}\eta, x)$,
for any vector field $x$ and for any $1$-forms $\xi$ and $\eta$.
Then there is a concise expression for the derivation $d_{P,\psi}$ 
which we now state.
\begin{proposition}
Let $d_P = [P, ~\cdot~]_{SN}$, where $\bra_{SN}$ is the Schouten-Nijenhuis
bracket.
The derivation $d_{P,\psi}$ on  $V^{\bullet}(M)$ is
\begin{equation}\label{diff}
d_{P,\psi}= d_P + i_{(\wedge^2 P^{\sharp}) (\psi)} \ . 
\end{equation}
\end{proposition}

\noindent{\it Proof}
Writing $P$ for $P^{\sharp}$, we compute, for any vector field $x$, and
any $1$-forms
$\xi$, $\eta$,
$$
(d_{P,\psi} x)(\xi, \eta) + \psi(P\xi,P\eta,x)
$$
$$
= L_{P\xi}<\eta,x> - L_{P\eta}<\xi,x> - <L_{P\xi}\eta,x> 
+ \, {\d}\xi(P\eta,x)
$$
$$
= <\eta,  L_{P\xi}x> - L_x <\xi, {P\eta}> - <\xi, L_{P\eta}x>
= - <\eta,  L_x{P\xi}> - <L_x \xi, {P\eta}> 
$$
$$
= - L_x<\eta, {P\xi}> + <L_x \eta, {P\xi}> - <L_x\xi, {P\eta}>
= - (L_xP)(\xi,\eta)=(d_P(x))(\xi,\eta).
$$
Moreover,
$i_{(\wedge^2 P^{\sharp}) (\psi)}$ is a derivation of $V^{\bullet}(M)$ 
which vanishes on $V^0(M)$ and coincides with $x \mapsto
\psi(P^{\sharp}~\cdot~,
P^{\sharp}~\cdot~,x)$ for $x \in V^1(M)$.
We have thus proved formula \eqref{diff}.

\medskip

Computing $[d_{P,\psi},d_{P,\psi}]$, we see that this derivation 
vanishes if and only
if condition \eqref{cns} is satisfied.
Thus, bracket $\bra^{P, \psi}$ is a Lie algebroid bracket if and only
if $P$ defines a Poisson structure with background $\psi$.
The operator $d_{P,\psi}$ is then 
the differential of the
Lie algebroid $(T^*M, \bra^{P,\psi})$. (See \cite{SW}.)

\medskip

In addition, the following morphism property of $P^{\sharp}$ is easily proved
\cite{yksESI}. 

\begin{proposition}
The relation
\begin{equation}
P^{\sharp}[\xi,\eta]^{P, \psi}= [P^{\sharp}\xi,P^{\sharp}\eta] \ .
\end{equation}
is equivalent to condition \eqref{cns}. 
\end{proposition}
\medskip

Further properties of the Poisson structures with background 
and their gauge
equi\-valence are studied in several recent articles, including
\cite{SW}, \cite{Rlmp}, \cite{yksESI} and
\cite{BC}.

\end{document}